\documentclass[article,12pt,a4paper,english,brazil]{article}

\usepackage[utf8]{inputenc}
\usepackage[T1]{fontenc}
\usepackage{indentfirst}
\usepackage{color}
\usepackage{graphicx}
\usepackage{microtype}
\usepackage{mathptmx}
\usepackage{etoolbox}
\usepackage{float}
\usepackage[bookmarks=false]{hyperref}
\usepackage{amsmath,amssymb,amsthm} 
\usepackage{graphicx}
\usepackage{wrapfig}
\usepackage{url}  

\title{Numerical convergence of a Telegraph Predator-Prey System}

\author{Kariston Stevan Luiz\thanks{Prof. Ms. Kariston Stevan Luiz, Depto. Matemática, UEL, Londrina, PR., Brasil, E-mail: kslpgmac@gmail.com}; 
Juniormar Organista\thanks{Prof. Ms. Juniormar Organista, Depto. Matemática, UEL, Londrina, PR., Brasil, E-mail: juniormarorganista@gmail.com}; 
Eliandro Rodrigues Cirilo\thanks{Prof. Dr. Eliandro Rodrigues Cirilo, Depto. Matemática, UEL, Londrina, PR., Brasil, E-mail: ercirilo@uel.br}; 
\\
\fontsize{13pt}{15.6pt} \selectfont
Neyva Maria Lopes Romeiro\thanks{Profa. Dra. Neyva Maria Lopes Romeiro, Depto. Matemática, UEL, Londrina, PR., Brasil, E-mail: nromeiro@uel.br}; 
Paulo Laerte Natti\thanks{Prof. Dr. Paulo Laerte Natti, Depto. Matemática, UEL, Londrina, PR., Brasil, E-mail: plnatti@uel.br} 
}

\begin{document}
	
\maketitle
	
Abstract: The numerical convergence of a Telegraph Predator-Prey system is studied. This system of partial differential equations (PDEs) can describe various biological systems with reactive, diffusive and delay effects. Initially, our problem is mathematically modeled. Then, the PDEs system is discretized using the Finite Difference method, obtaining a system of equations in the explicit form in time and implicit form in space. The consistency of the Telegraph Predator-Prey system discretization was verified. Next, the von Neumann stability conditions were calculated for a Predator-Prey system with reactive terms and for a Telegraph system with delay. For our Telegraph Predator-Prey system, through numerical experiments, it was verified tat the mesh refinement and the model parameters (reactive constants, diffusion coefficient and delay term) determine the stability/instability conditions of the model.
	
	\vspace{0.15cm}

\newpage

\section*{Introduction} 

Currently, there is a growing interest in the study of population dynamics, mainly due to the need to have better control of epidemics, or to reduce the economic, biological and social losses caused by invasive species, among other needs.

The first mathematical studies aimed at describing interactions between populations took place in the 1920s \cite{intro}. During this period, the Lotka-Volterra mathematical model emerged \cite{lotka, volterra}. This model describes the predator-prey interaction between two populations through the system of ordinary differential equations,
\begin{eqnarray}
\frac{dS_1}{dt} &=& a_1S_1 - c_1S_1S_2 \nonumber \\ \\
\frac{dS_2}{dt} &=& -a_2S_2 + c_2S_1S_2, \nonumber
\label{1}
\end{eqnarray}
\noindent
where $S_1 > 0$ and $S_2 > 0$ denote the population densities of the interacting species, with $S_1$ being the density of prey and $S_2$ the density of predators. Furthermore, $a_1 > 0$ is the birth rate of the species $S_1$, $a_2 > 0$ is the death rate of the species $S_2$, and the parameters $c_1 > 0$ and $c_2 > 0$ are the interaction rates between the two species. The non-derivative terms on the right side of (1) are called reactive terms.

Later in the 1950s, C. S. Holling carried out experiments to investigate how the rate of prey capture by a predator is related to the density of the prey \cite{holling1,holling2}, a relationship called the functional response.

\begin{figure}[ht]
    \begin{center}
    \caption{The three types of functional responses identified by Holling.}
    \includegraphics[width=8cm]{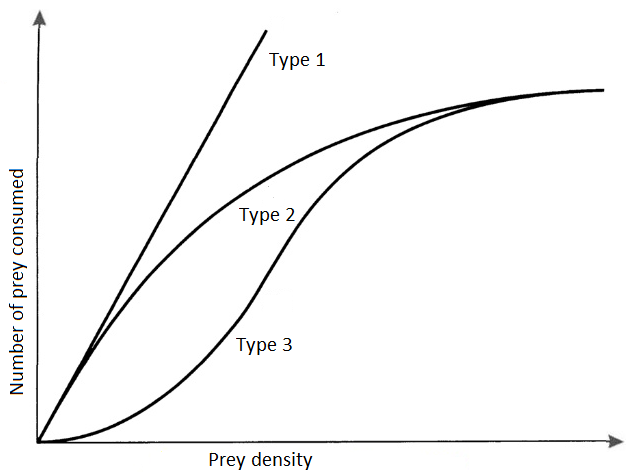}
    \label{holl}
    \end{center}
\end{figure}

Holling identified three general categories of functional responses (Figure \ref{holl}). The type 1 functional response is linear, when the number of prey consumed increases in direct proportion to prey density. The type 2 functional response says that as the prey population increases, predators become satiated and consume a constant number of prey (saturation). The type 3 functional response is similar to type 2, except at low prey density, when prey switching effects occur \cite{Tansky}.

Still in the 1950s, Carl B. Huffaker investigated the effects of spatial structure on the interaction of mite populations. Among several results, Huffaker showed that the predator-prey interaction could not survive in a homogeneous environment without dispersion \cite{huff1}. 

Currently, through more complex mathematical modeling, there are works that emphasize the studies of the dynamics of invasive species, epidemics and other biological phenomena \cite{Deroos,pp3,pp1,petro1,pp4,cirilo1}.

As the description of the effects of predation on one species by another improved, the mathematical models became more complex. Thus, in parallel with the development of more realistic and complex biological models, studies on the convergence of numerical methods applied to these models have become essential. 

The Lax Equivalence Theorem is fundamental for the analysis of the convergence of numerical solutions of PDEs. The theorem says that in a well-posed initial value problem that has been consistently discretized, stability of the numerical method is a necessary and sufficient condition for numerical convergence \cite{lalax}. 

About stability of numerical methods, it is associated with the propagation of numerical errors in the iterative process \cite{hhh}. The numerical method is said to be stable if the errors decrease along the iterative process, otherwise, if the errors increase, the numerical solution will diverge and the numerical method is said to be unstable. Finally, a numerical method is said to be conditionally stable if it depends on certain parameters so that the errors remain limited. In the case of linear problems, the Von Neumann stability analysis is widely used  \cite{cuminato1999discretizaccao}. In the case of nonlinear initial condition problems, the stability of a numerical method can hardly be verified analytically \cite{est1}. In such situations, numerical experiments are used to analyze the stability of the numerical method used \cite{natti2}.

About consistency of numerical method, we say that a discretized problem is consistent when it tends to the original differential equation if the increments in the independent variables tend to zero \cite{consite1}.

In this context, the objective of this work is to verify the consistency and stability of a numerical method applied to a Telegraph Predator-Prey system. The conditions for the convergence of the numerical model are analyzed through numerical experiments.

\section*{Modeling the Telegraph Predator-Prey system}

In this section, the Maxwell-Cattaneo Diffusion equation system is first modeled and then the Telegraph Predator-Prey equation system is developed.

\subsection*{Maxwell-Cattaneo Diffusion system}

The classical law of diffusion, Fick's law \cite{paul2014fick}, is a law that describes the diffusion of a property in systems that are not in equilibrium. In systems where there are concentration gradients of a property, then there is a flow of this property which tends to homogenize the system.This homogenizing flow will go in the opposite direction of the gradient and, if this flow is weak, it can be approximated by the first term of a Taylor series, resulting in Fick's law. In the case of a Predator-Prey system in one dimension, this modeling is expressed mathematically by \cite{paul2014fick},
\begin{equation}\label{equ1leidefick}
J_{j} = -D_{j} \frac{\partial S_{j}}{\partial x},
\end{equation}
where flow $J_1$ is due to diffusion of prey $S_1$ and flow $J_2$ is due to diffusion of predators $S_2$, while $D_1$ and $D_2$ are the diffusivity coefficients for prey and predators, respectively. Note that the negative sign in \ref{equ1leidefick} indicates that the homogenizing flow occurs in the opposite direction of the concentration gradient, ie, from high concentrations to low concentrations.

On the other hand, Fick's law proposes that signals propagate with infinite speed, which in practice does not happen, configuring the so-called "Paradox of Fourier's law" \cite{mickens}. To correct this problem, Fick's law is modified by the Maxwell-Cattaneo Diffusion law \cite{cata} which introduces a lag time $\tau$ for each action. Thus, due to the finite speed of propagation of information, the system does not react instantly to an action.
 
For predator-prey systems, the Maxwell-Cattaneo Diffusion law is given by \cite{cirilo1}
 \begin{equation}
	\displaystyle \left(1+ \tau_{j} \frac{\partial }{\partial t}\right)J_{j} = -D_{j} \frac{\partial S_{j}}{\partial x},
	\label{EQ:M-C1}
      \end{equation}
where $\tau_1$ is the reaction time of the prey when exposed to predation and $\tau_2$ is the predator's reaction time to capture prey. Note that, contrary to the equation system (1), in the system (3) we have $S_j = S_j(x,t)$, for $j=1.2$, so that now the population densities of prey and predator, respectively, also depend on the spatial variable $x$.

\subsection*{Telegraph Predator-Prey system}

For predator-prey systems, in the context of Huffaker's hypotheses, the Maxwell-Cattaneo diffusion law must be incorporated into the population density conservation law. According to the Conservation Principle, the rate of change of a property in a volume $V$ must be equal to the net flux of that property through the surface of $V$, plus the amount of the property transformed into the interior of V due to reactive effects \cite{de2000tecnicas}. Thus, the conservation equation for predator-prey system, considering diffusive, reactive and delay effects, is written as
 \begin{equation}
	\displaystyle \frac{\partial S_j}{\partial t} = -\frac{\partial J_j}{\partial x} + F_j\left(S_j\right),
	\label{EQ:R}
      \end{equation}
where $J_j$ are the population flows described by (\ref{EQ:M-C1}) and $F_j(t,x,S)$ represents the prey reaction term, if $j=1$, and the predator reaction term, if $j=2$. Note that the equation (\ref{EQ:R}) is a transport equation.

Therefore, deriving \ref{EQ:R} with respect to time $t$ and \ref{EQ:M-C1} with respect to coordinate $x$, we have
\begin{equation}
	\displaystyle \frac{\partial^{2} S_j}{\partial t^{2}} = -\frac{\partial}{\partial t}\left(\frac{\partial J_j}{\partial x}\right) + \frac{\partial}{\partial t}F_j\left(S_j\right)
	\label{EQ:R-D2}
      \end{equation}
      
\begin{equation}
	\displaystyle \tau_j\frac{\partial}{\partial x}\left(\frac{\partial J_j}{\partial t}\right) +\frac{\partial J_j}{\partial x} = -D_j\frac{\partial^{2} S_j}{\partial x^{2}}.
	\label{EQ:M-C2}
      \end{equation}
\noindent      
Noting that $F_j = F_j(S_j,x,t)$ and $S_j = S_j(x,t)$, it follows from the chain rule that
\begin{equation}
	\displaystyle \frac{\partial}{\partial t}F_j\left(S_j\right) = \frac{d}{d S} F_j(S_j)\hspace{1.5mm} \frac{\partial S_j}{\partial t},
	\label{EQ:R-DO2}
      \end{equation}
so that \ref{EQ:R-D2} is rewritten as
\begin{equation}
	\displaystyle \frac{\partial^{2} S_j}{\partial t^{2}} = -\frac{\partial}{\partial t}\left(\frac{\partial J_j}{\partial x}\right) + \frac{d}{d S_j} F_j(S_j)\hspace{1.5mm}\frac{\partial S_j}{\partial t}.
	\label{EQ:R-D3}
      \end{equation}
\noindent     
Multiplying the equation \ref{EQ:R-D3} by $\tau_j$ and subtracting it from the equation \ref{EQ:M-C2}, we get that
 \begin{equation}
	\displaystyle \tau_j\frac{\partial^2 S_j}{\partial t^2} -
	              \tau_j \frac{d}{d S_j}F_j\left(S_j\right)\hspace{1.5mm} \frac{\partial S_j}{\partial t} =
	              D_j\frac{\partial^2 S_j}{\partial x^2} + \frac{\partial J_j}{\partial x}.
	\label{EQ:RT}
      \end{equation}
Finally, obtaining $\frac{\partial J_j}{\partial x}$ from \ref{EQ:R} and substituting in \ref{EQ:RT}, we obtain the delayed predator-prey equations, whose population densities are subject to reactive-diffusive processes \cite{equacao1,cirilo1},i.e.,
  \begin{equation}
	\displaystyle \tau_j\frac{\partial^2 S_j}{\partial t^2} + 
	              \left[ 1 - \tau_j \frac{d}{d S_j}F_j\left(S_j\right)\right]\frac{\partial S_j}{\partial t} =
	              D_j\frac{\partial^2 S_j}{\partial x^2} + F_j\left(S_j\right).
	\label{EQ:R-T}
      \end{equation}
The equations \ref{EQ:R-T} will be designated in the remaining of the work as Telegraph Predator-Prey equations. Note that the equations \ref{EQ:R-T} have the same structure as the Telegraph equation, which was derived by William Thomson (Lord Kelvin) to describe the propagation of electrical signals in long conducting cables.

\section*{Discretization of the Telegraph Predator-Prey system}

Consider the equations \ref{EQ:R-T} with the following initial and boundary conditions
\begin{equation}
	\tau_{1}\frac{\partial^2 S_{1}}{\partial t^2}  + 
	\left[ 1 - \tau_{1}\frac{d F_{1}}{d S_{1}}\right]\frac{\partial S_{1}}{\partial t} =  D_{1}\frac{\partial^2 S_{1}}{\partial x^2} + F_{1}\label{EQ:R_T_C}
\end{equation}
\begin{equation}
	\tau_{2}\frac{\partial^2 S_{2}}{\partial t^2}  + 
	\left[ 1 - \tau_{2} \frac{d F_{2}}{d S_{2}}\right]\frac{\partial S_{2}}{\partial t} =  D_{2}\frac{\partial^2 S_{2}}{\partial x^2} + F_{2}\label{EQ:R_T_C_1} 
\end{equation}
\begin{equation}
	S_{j}(x,0)=S_{j}^{0} \; \; , \; \; \frac{\partial S_{j}(x,t)}{\partial t}\bigg|^{t=0} = 0 \; \; , \; \; \forall x\in[0,L] \label{condicaoinicial}
\end{equation}
\begin{equation}
	S_{j}(0,t)=S_{j}(L,t) = 0 \; \; , \; \; \forall t\in[0,T] \; \; , \; \; j=1,2, \label{condicaodecontorno}
\end{equation}
where $t$ and $x$ are the temporal and spatial variables, $\tau_{1}$ and $\tau_{2}$ are the delay parameters of the populations, $D_{1}$ and $D_{2} $ are the diffusivity coefficients of the populations, $S_{1}(x,t)$ and $S_{2}(x,t)$ are the population densities, finally $F_{1}$ and $F_{2} $ are the prey and predator reactive terms of the populations, respectively. 
In this work, it is considered that the reactive terms present functional responses of the Holling 1 type (see Figure 1), i.e.,
\begin{equation}\label{termoreativopresa}
	F_{1} = F_{1}(S_{1},S_{2}) = a_{1}S_{1} - b_{1}S_{1}^{2} - c_{1}S_{1}S_{2}
\end{equation}
\begin{equation}\label{termoreativopredador}
	F_{2} = F_{2}(S_{1},S_{2}) = - a_{2}S_{2} + c_{2}S_{1}S_{2},
\end{equation}
where $a_{1}$ is the prey birth rate, $b_{1}$ is the prey saturation term, $c_{1}$ is the prey death rate due to predation, $a_{2 }$ is the predator death rate in the absence of prey and $c_{2}$ is the predator birth rate due to predation \cite{natti1}. 

The initial conditions $S_{1}^{0}$ and $S_{2}^{0}$ are the initial population densities of prey and predators, respectively. Finally, the  Dirichlet boundary conditions (14) impose that at the boundary of the problem domain, the population densities are null.

\begin{figure*}[ht]
	\begin{center}
	\caption{Discrete one-dimensional mesh used for the 
Telegraph Predator-Prey system.}
		\includegraphics[width=13cm]{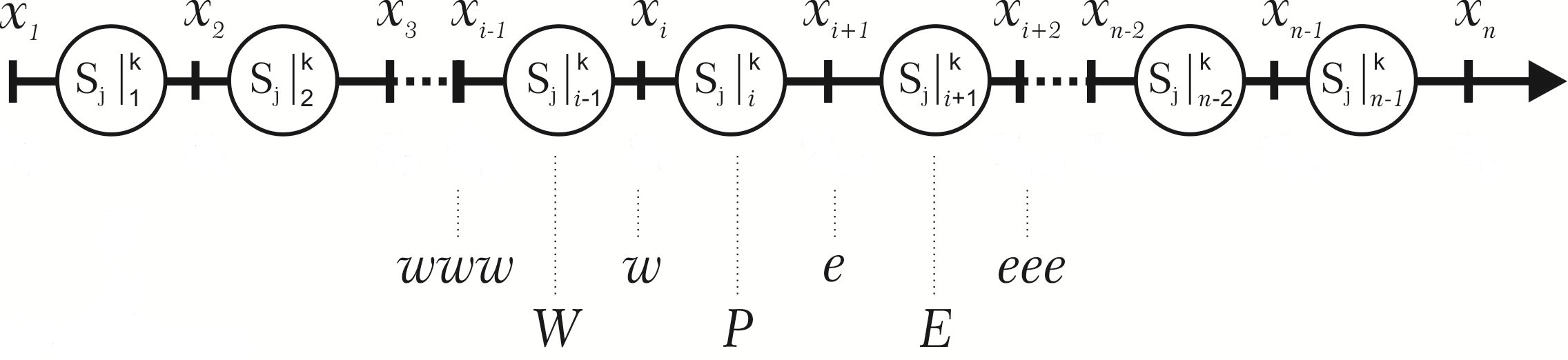}
\label{deslocada} \vspace{4pt}
\end{center}
\end{figure*}

To discretize the Telegraph Predator-Prey equations (11-
12) the one-dimensional mesh shown in Figure \ref{deslocada} is used.

Cardinal point notation is used in the discretization of the problem domain. The labels $P$, $W$ and $E$ stand for center ($P$ point where the calculation is performed), east and west, respectively. Lowercase acronyms are cardinal variations from the center $P$. As shown in Figure \ref{deslocada}, the population densities $S_1$ and $S_2$ are calculated at the center of the cell. Note that storage shifted to population densities has a positive impact on numerical computation, reducing numerical instability \cite{paolo1,de2000tecnicas,paolo2,barba1,tati1}.

As for the numerical method for discretizing the PDEs, the Finite Difference Technique was used  \cite{cuminato1999discretizaccao,burden2008analise}. Thus, central finite differences are used for the discretization of the second derivatives, that is,
\begin{eqnarray}
\left.\left(\frac{\partial^{2}S_{j}}{\partial t^{2}}\right)\right|_{P}^{k}  \approx \frac{1}{(\Delta t)^{2}}\left(S_{j}|_{P}^{k + 1} - 2S_{j}|_{P}^{k} + S_{j}|_{P}^{k - 1}\right)
\label{eqsjtt} 
\end{eqnarray}
\vspace{-0.4cm}
\begin{equation}
\label{eqsjxx}
\left.\left(\frac{\partial^{2}S_{j}}{\partial x^{2}}\right)\right|_{P}^{k} \approx \frac{1}{(\Delta x)^{2}}\left(\left. S_{j}\right|_{E}^{k} - 2\left.S_{j}\right|_{P}^{k} + \left.S_{j}\right|_{W}^{k} \right),
\end{equation}
while for the first derivative in time, regressive finite differences are used
\begin{equation}\label{eqsjt}
\left.\left(\frac{\partial}{\partial t}S_{j}\right)\right|_{P}^{k} \approx \frac{1}{\Delta t}\left(\left.S_{j}\right|_{P}^{k} - \left.S_{j}\right|_{P}^{k-1}\right).
\end{equation}

Finally, for the equation \ref{EQ:R_T_C}, it is still necessary to discretize the term 
\begin{eqnarray}
\left.\left(\left[ 1 - \tau_{1} \frac{\partial}{\partial S_{1}}F_{1}\left(S_{1},S_{2}\right)\right]\frac{\partial S_{1}}{\partial t}\right)\right|_{P}^{k} =&& \nonumber \\
\left.\left(\frac{\partial S_{1}}{\partial t}\right)\right|_{P}^{k} - \tau_{1}\left.\left(\frac{\partial}{\partial S_{1}}F_{1}\left(S_{1},S_{2}\right)\right)\right|_{P}^{k}\left.\left(\frac{\partial S_{1}}{\partial t}\right)\right|_{P}^{k} . && 
 \label{eqs1f1t}
\end{eqnarray}
From \ref{termoreativopresa}, it follows that
\begin{equation}\label{eqs1f1t_2}
\left.\left(\frac{\partial}{\partial S_{1}}F_{1}\left(S_{1},S_{2}\right)\right)\right|_{P}^{k} = a_{1} - 2b_{1}\left.S_{1}\right|_{P}^{k} - c_{1}\left.S_{2}\right|_{P}^{k}, \end{equation}
then using the equations \ref{eqsjt} and \ref{eqs1f1t_2} in the equation \ref{eqs1f1t}, one obtains that
\begin{eqnarray} 
\left.\left(\frac{\partial}{\partial S_{1}}F_{1}\left(S_{1},S_{2}\right)\frac{\partial S_{1}}{\partial t}\right)\right|_{P}^{k}  \approx 
\nonumber \\
\left(a_{1} - 2b_{1}\left.S_{1}\right|_{P}^{k} - c_{1}\left.S_{2}\right|_{P}^{k}\right)\left(\frac{1}{\Delta t}\left(\left.S_{1}\right|_{P}^{k} - \left.S_{1}\right|_{P}^{k-1}\right)\right) 
\approx \nonumber\\
\frac{1}{\Delta t}\left(a_{1}\left.S_{1}\right|_{P}^{k} - 2b_{1}\left.S_{1}^{2}\right|_{P}^{k} - c_{1}\left.S_{1}\right|_{P}^{k}\left.S_{2}\right|_{P}^{k}\right.\nonumber\\
- \left. a_{1}\left.S_{1}\right|_{P}^{k-1} + 2b_{1}\left.S_{1}\right|_{P}^{k-1}\left.S_{1}\right|_{P}^{k} + c_{1}\left.S_{1}\right|_{P}^{k-1}\left.S_{2}\right|_{P}^{k}\right).
\label{eqs1f1t_3}
\end{eqnarray}
Substituting (15), \ref{eqsjtt}, \ref{eqsjxx}, \ref{eqsjt}, and \ref{eqs1f1t_3} in \ref{EQ:R_T_C}, it follows that the discretization of the prey equation is given by
\begin{eqnarray}
 S_1|_{P}^{k + 1}  &=&  \Omega_{1}\left( \Upsilon_{1}S_{1}|_{W}^{k} + \Pi_{1}S_{1}^{2}|_{P}^{k} \right.\nonumber \\
 &+& \left.\Phi_{1}S_{1}|_{P}^{k} + \Lambda_{1}S_{1}|_{E}^{k} + \Gamma_{1}\right),
\label{finalpresa}
\end{eqnarray}
where
\begin{eqnarray}
\Upsilon_{1}  = - \frac{D_1}{{\Delta x}^2} \;\;\;\;,\;\;\;\;
\Lambda_{1} = - \frac{D_1}{{\Delta x}^2}
\end{eqnarray}
\begin{eqnarray}
\Omega_1 = - \frac{{\Delta t}^2}{\tau_1} \;\;\;\;,\;\;\;\;
\Pi_{1}  =  \frac{2\tau_1b_1}{\Delta t} + b_{1}
\end{eqnarray}
\begin{eqnarray}
\Phi_{1} & = &  \frac{1 - \tau_1\left(a_1 - c_1S_{2}|_{P}^{k} +2b_1S_{1}|_{P}^{k-1}\right)}{\Delta t}\nonumber\\
&-& \frac{2\tau_1}{{\Delta t}^2}
+ \frac{2D_1}{{\Delta x}^2} - a_{1} + c_{1}S_{2}|_{P}^{k}
\end{eqnarray}
\begin{eqnarray}
\Gamma_{1}=\frac{\tau_1S_1|_{P}^{k - 1}}{{\Delta t}^2} - \frac{\left(1 - \tau_1\left.(a_{1} - c_{1}S_{2}|_{P}^{k})\right.\right)S_{1}|_{P}^{k - 1}}{\Delta t}.
\end{eqnarray}
Note that the discretized prey equation \ref{finalpresa} is explicit in time and implicit in space.

In a similar way, the formula for the discretization of the predator equation is obtained. The only difference is in the derivation of the predator reactive term. From equation (16)
\begin{equation}\label{eqs1f2t_2}
\left.\left(\frac{\partial}{\partial S_{2}}F_{2}\left(S_{1},S_{2}\right)\right)\right|_{P}^{k} = -a_{2} + c_{2}\left.S_{1}\right|_{P}^{k}. 
\end{equation}
So the equation for discretizing the predator equation is given by
\begin{eqnarray}
S_2|_{P}^{k + 1}  =  \Omega_{2}\left( \Upsilon_{2}S_{2}|_{W}^{k}  +\Phi_{2}S_{2}|_{P}^{k} + \Lambda_{2}S_{2}|_{E}^{k} + \Gamma_{2}\right),\label{finalpredador}
\end{eqnarray}
where
\begin{eqnarray}
\Upsilon_{2} = - \frac{D_2}{{\Delta x}^2} \;\;\;\;,\;\;\;\;
\Lambda_{2} = - \frac{D_2}{{\Delta x}^2} \;\;\;\;,\;\;\;\;
\Omega_2 = - \frac{{\Delta t}^2}{\tau_2}
\end{eqnarray}
\begin{eqnarray}
\Phi_{2} & = & \frac{1 - \tau_2(-a_2 + c_2S_1|_{P}^{k})}{\Delta t}
\nonumber\\
&-& \frac{2\tau_2}{{\Delta t}^2} + \frac{2D_2}{{\Delta x}^2}
+ a_2 - c_2S_1|_{P}^{k}
\end{eqnarray}
\begin{eqnarray}
\Gamma_{2}=\frac{\tau_2S_2|_{P}^{k - 1}}{{\Delta t}^2} - \frac{\left(1 - \tau_2(-a_2 +c_2S_1|_{P}^{k})\right)S_2|_{P}^{k -1}}{\Delta t}.
\end{eqnarray}
Again, note that the discretized prey equation \ref{finalpredador} is explicit in time and implicit in space.

\section*{Numerical model consistency}

The numerical solution of a problem does not always tend to the exact solution of the same problem. The resulting error between the exact solution of the continuous problem and the numerical solution obtained from the discretized equations is called Local Truncation Error (LTE) \cite{cuminato1999discretizaccao}. In most of the applied problems, the exact solution is not known and this is our case. In such situations one can estimate the Local Truncation Error through Taylor Series \cite{burden2008analise} and use this estimate to prove that the method is consistent at the limit of the continuum, when $\Delta x,\Delta t\to 0$. 

Consider EDPs \ref{EQ:R_T_C} and \ref{EQ:R_T_C_1}. Using the Finite Difference method, we have the following discretized equations
\begin{eqnarray}
\frac{\tau_1}{\Delta t^{2}}\left(S_{1}\bigg|_{P}^{k + 1}  - 2S_{1}\bigg|_{P}^{k} + S_{1}\bigg|_{P}^{k - 1}\right)  
\;\;\;\;\;\;\;\;\;\;\;\;\;
\nonumber\\
+ \frac{1}{\Delta t}\left[\left(1 - \tau_{1}\left(a_{1} - c_{1}\left.S_{2}\right|_{P}^{k} - 2b_{1}\left.S_{1}\right|_{P}^{k}\right)\right)S_{1}|_{P}^{k}\right.
\nonumber\\
- \left.\left(1 - \tau_{1}\left(a_{1} - c_{1}\left.S_{2}\right|_{P}^{k} - 2b_{1}\left.S_{1}\right|_{P}^{k}\right)\right)\left.S_{1}\right|_{P}^{k-1}\right] 
\nonumber\\
- \frac{D_1}{\Delta x^{2}}\left( S_{1}\bigg|_{E}^{k} - 2S_{1}\bigg|_{P}^{k} + S_{1}\bigg|_{W}^{k} \right) - F_1\bigg|_P^{k} = 0
\label{final1}
\end{eqnarray}
and
\begin{eqnarray}
\frac{\tau_2}{\Delta t^{2}}\left(S_{2}\bigg|_{P}^{k + 1} - 2S_{2}\bigg|_{P}^{k} + S_{2}\bigg|_{P}^{k - 1}\right) 
\;\;\;\;\;\;\;\;\;\;\;
\nonumber\\
+ \frac{1}{\Delta t}\left[\left(1 - \tau_{2}\left(- a_{2} + c_{2}\left.S_{1}\right|_{P}^{k}\right)\right)S_{2}|_{P}^{k}\right.
\;\;\;\;\;\;
\nonumber\\
- \left.\left(1 - \tau_{2}\left(- a_{2} + c_{2}\left.S_{1}\right|_{P}^{k}\right)\right)S_{2}|_{P}^{k-1}\right] 
\;\;\;\;\;\;
\nonumber\\
- \frac{D_2}{\Delta x^{2}}\left( S_{2}\bigg|_{E}^{k} - 2S_{2}\bigg|_{P}^{k} + S_{2}\bigg|_{W}^{k} \right) - F_2\bigg|_P^{k} = 0,
\label{final2}
\end{eqnarray}
where $F_1$ and $F_2$ are the prey and predator reactive terms, equations (15) and (16), respectively.

Assuming that $W = P - \Delta x$ and $E = P + \Delta x$, let the Taylor Series expansions of the population densities in (33) and (34) be calculated around $k$ and $P $, that is,
\begin{equation}
\label{analise_consistencia_s1_termo_final1} 
	S_{j}\bigg|_{P}^{k + 1} = S_{j}\bigg|_{P}^{k} + \Delta t\frac{\partial S_{j}}{\partial t}\bigg|_{P}^{k} + \frac{(\Delta t)^{2}}{2}\frac{\partial^{2} S_{j}}{\partial t^{2}}\bigg|_{P}^{k} + \frac{(\Delta t)^{3}}{3!}\frac{\partial^{3} S_{j}}{\partial t^{3}}\bigg|_{P}^{k} + \mathcal{O}(\Delta t^{4})
\end{equation} 
\begin{equation}
\label{analise_consistencia_s1_termo_final2}
	S_{j}\bigg|_{P}^{k - 1} = S_{j}\bigg|_{P}^{k} - \Delta t\frac{\partial S_{j}}{\partial t}\bigg|_{P}^{k} + \frac{(\Delta t)^{2}}{2}\frac{\partial^{2} S_{j}}{\partial t^{2}}\bigg|_{P}^{k} - \frac{(\Delta t)^{3}}{3!}\frac{\partial^{3} S_{j}}{\partial t^{3}}\bigg|_{P}^{k} + \mathcal{O}(\Delta t^{4})
\end{equation} 
\begin{equation}
\label{analise_consistencia_s1_termo_final3}
	S_{j}\bigg|_{W}^{k} = S_{j}\bigg|_{P}^{k} - \Delta x\frac{\partial S_{j}}{\partial x}\bigg|_{P}^{k} + \frac{(\Delta x)^{2}}{2}\frac{\partial^{2} S_{j}}{\partial x^{2}}\bigg|_{P}^{k} - \frac{(\Delta x)^{3}}{3!}\frac{\partial^{3} S_{j}}{\partial x^{3}}\bigg|_{P}^{k} + \mathcal{O}(\Delta x^{4})
\end{equation} 
\begin{equation}
\label{analise_consistencia_s1_termo_final4}
	S_{j}\bigg|_{E}^{k} = S_{j}\bigg|_{P}^{k} + \Delta x\frac{\partial S_{j}}{\partial x}\bigg|_{P}^{k} + \frac{(\Delta x)^{2}}{2}\frac{\partial^{2} S_{j}}{\partial x^{2}}\bigg|_{P}^{k} + \frac{(\Delta x)^{3}}{3!}\frac{\partial^{3} S_{j}}{\partial x^{3}}\bigg|_{P}^{k} + \mathcal{O}(\Delta x^{4})
\end{equation}

\noindent Substituting the equations \ref{analise_consistencia_s1_termo_final1}- \ref{analise_consistencia_s1_termo_final4} in the equations \ref{final1} and \ref{final2}, with their respective $j = 1, 2$, we obtain that
\begin{eqnarray}
\underbrace{\tau_1\frac{\partial^{2} S_{1}}{\partial t^{2}}\bigg|_{P}^{k} + \left[1 - \tau_1	\frac{d F_{1}}{d S_{1}}\bigg|_{P}^{k}\right]\frac{\partial S_1}{\partial t}\bigg|_P^{k} - D_1\frac{\partial^2S_1}{\partial x^2}\bigg|_P^{k} - F_1\bigg|_P^{k}}_{\mbox{PDE}}
\nonumber\\
= \underbrace{\frac{\Delta t}{2}\left(1 - \tau_1\frac{d F_{1}}{d S_{1}}\bigg|_{P}^{k}\right)\frac{\partial^2S_1}{\partial t^2}\bigg|_P^{k}
  -  \frac{\Delta t^2}{3!}\left(1 - \tau_1\frac{d F_{1}}{d S_{1}}\bigg|_{P}^{k}\right)\frac{\partial^3S_1}{\partial t^3}\bigg|_P^{k}}_{\mbox{Local Truncation Error}}
  \nonumber\\
+ \underbrace{ \left(1  -\tau_1\frac{d F_{1}}{d S_{1}}\bigg|_{P}^{k}\right)\mathcal{O}(\Delta t^{3})
 -   2\tau_1\mathcal{O}(\Delta t^{2}) + 2D_1\mathcal{O}(\Delta x^{2})}_{\mbox{Local Truncation Error}}
 \label{vv1_1}
\end{eqnarray}
and
\begin{eqnarray}
\underbrace{\tau_2\frac{\partial^{2} S_{2}}{\partial t^{2}}\bigg|_{P}^{k} + \left[1 - \tau_2	\frac{d F_{2}}{d S_{2}}\bigg|_{P}^{k}\right]\frac{\partial S_2}{\partial t}\bigg|_P^{k} - D_2\frac{\partial^2S_2}{\partial x^2}\bigg|_P^{k} - F_2\bigg|_P^{k}}_{\mbox{PDE}}
\nonumber\\
 = \underbrace{\frac{\Delta t}{2}\left(1 - \tau_2\frac{d F_{2}}{d S_{2}}\bigg|_{P}^{k}\right)\frac{\partial^2S_2}{\partial t^2}\bigg|_P^{k}
  -  \frac{\Delta t^2}{3!}\left(1 - \tau_2\frac{d F_{2}}{d S_{2}}\bigg|_{P}^{k}\right)\frac{\partial^3S_2}{\partial t^3}\bigg|_P^{k}}_{\mbox{Local Truncation Error}}
  \nonumber\\
+ \underbrace{ \left(1  -\tau_2\frac{d F_{2}}{d S_{2}}\bigg|_{P}^{k}\right)\mathcal{O}(\Delta t^{3})
 -   2\tau_2\mathcal{O}(\Delta t^{2}) + 2D_2\mathcal{O}(\Delta x^{2}).}_{\mbox{Local Truncation Error}}
 \label{vv1_2}
\end{eqnarray}
Finally, taking the limit of the continuum, when $\Delta t,\Delta x\to 0$, the equations \ref{vv1_1} and \ref{vv1_2} tend to
\begin{equation}
\label{consistentelimite1}
    \tau_1\frac{\partial^{2} S_{1}}{\partial t^{2}}\bigg|_{P}^{k} + \left[1 - \tau_1	\frac{d F_{1}}{d S_{1}}\bigg|_{P}^{k}\right]\frac{\partial S_1}{\partial t}\bigg|_P^{k} - D_1\frac{\partial^2S_1}{\partial x^2}\bigg|_P^{k} - F_1\bigg|_P^{k} = 0
\end{equation}
and 
\begin{equation}
\label{consistentelimite2}
    \tau_2\frac{\partial^{2} S_{2}}{\partial t^{2}}\bigg|_{P}^{k} + \left[1 - \tau_2	\frac{d F_{2}}{d S_{2}}\bigg|_{P}^{k}\right]\frac{\partial S_2}{\partial t}\bigg|_P^{k} - D_2\frac{\partial^2S_2}{\partial x^2}\bigg|_P^{k} - F_2\bigg|_P^{k} = 0.
\end{equation}
Note that the equations (\ref{consistentelimite1}) and  (\ref{consistentelimite2}) are the PDEs (11) and (12), calculated at the mesh point $P$ and at the instant $k$. It is concluded that the discretized equations (23) and (29), or equivalently equations \ref{final1} and \ref{final2}, are consistent with the PDEs (11) and (12), respectively.

\section*{Numerical model stability} 
The Von Neumann stability condition is based on the superposition principle, that is, which the error is the superposition of the errors accumulated at each iteration \cite{vneumann1}. The Von Neumann stability condition produces a necessary but not sufficient condition for stability \cite{cuminato1999discretizaccao}. 

Consider $\Delta x$ and $\Delta t$ the partitions in space and time, respectively, and $I = \sqrt{-1}$ the imaginary complex number. Let $E_i$, for $i = 1, ... N$, be the  error at each mesh point in time step $t = 0$. Then we write $E_i$ through a complex Fourier series \cite{de2000tecnicas}
\begin{equation}
E_i = \sum_{n=1}^{N}a_{n}^{0} \;\; e^{I\alpha_{n}
(i\Delta x)}
\label{vn88}, \hspace{7mm} \hbox{i = 1, ... N,}
\end{equation}
where $\alpha_{n} = \frac{n\pi}{L}$ and $N\Delta x = L$.

The time dependence of the error $E_i$ is incorporated through the amplitude $a_n^0 \rightarrow a_n^k$. It is assumed that the error tends to grow or decay exponentially with time, so the error \ref{vn88} in time $k$ is written as
\begin{equation}
    E_i^k = \sum_{n=0}^{N}e^{\gamma_n k}e^{I\alpha_{n}i\Delta x}.\label{vn99}
\end{equation}
Now $E_i^k$ is the error at the position $i$ and at the time step $k$.

The error equation \ref{vn99} is a propagation of harmonic waves of the type $e^{\gamma k}e^{I\xi i}$. In order for this propagation to be stable, the absolute value of the error amplitude must be less than or equal to unity \cite{cuminato1999discretizaccao, de2000tecnicas}, i.e.,
\begin{eqnarray}
 |e^{\gamma k}| \leq 1.\label{pratico}
 \end{eqnarray}

Next, the Von Neumann stability conditions will be obtained for the equations \ref{EQ:R_T_C} - \ref{EQ:R_T_C_1} when:

\vspace{0,2cm}
\noindent
Case I) \hspace{0,3cm} $F_1 , F_2 \neq 0$ and $\tau_{1} , \tau_{2} = 0$ (Predator-Prey system with reactive terms).

\vspace{0,2cm}
\noindent
Case II) \hspace{0,3cm} $F_1 , F_2 = 0$, $\tau_{1} , \tau_{2} \neq 0$  (Telegraph equations). 

\vspace{0,2cm}
\noindent
Case III) \hspace{0,3cm} $F_1 , F_2 \neq 0$ and $\tau_{1} ,  \tau_{2}\neq 0$ (Telegraph Predator-Prey system). In this case the stability analysis will be obtained through numerical experiments.

\subsection*{Von Neumann conditions for a Predator-Prey system}\label{sectioncalor2}
Consider the following Diffusive Predator-Prey system
\begin{eqnarray}
	\frac{\partial S_{1}}{\partial t}  = D_1\frac{\partial^2 S_{1}}{\partial x^2} + a_1S_1 - b_1S_1^2 - c_1S_1S_2
	\label{calor1} 
	\\
	\frac{\partial S_{2}}{\partial t}  = D_2\frac{\partial^2 S_{2}}{\partial x^2} - a_2S_2 +  c_2S_1S_2,\;\;\;\;\;\;
	\label{calor2}
\end{eqnarray}
where $S_1$, $S_2$ are the densities of prey and predator populations, respectively. The parameters $a_1$, $a_2$, $c_1$, $c_2$, and $b_1$ are positive constants described in (15) and (16), while $D_1>0$ and $D_2>0$ are the prey and predator diffusion rates, respectively.  Note that the equations \ref{calor1} and \ref{calor2} are the equations \ref{EQ:R_T_C} and \ref{EQ:R_T_C_1} when $\tau_1 = \tau_2 = 0$.

For the discretization of the equations \ref{calor1} and \ref{calor2}, we use progressive finite differences in time first-order derivative and central finite differences in space second-order derivative. Then the discretized versions of the equations \ref{calor1} and \ref{calor2} have the form
\begin{eqnarray}
S_1|_P^{k+1} = S_1|_{P}^{k} + \sigma_1\left(S_1|_E^{k} - 2S_1|_P^{k} + S_1|_W^{k}\right) \nonumber \\
+ \Delta tS_1|_P^{k}\left(a_1 - b_1S_1|_{P}^{k} - c_1S_2|_{P}^{k}\right) \label{calor1dis}
\\
S_2|_P^{k+1} = S_2|_{P}^{k} + \sigma_2\left(S_2|_E^{k} - 2S_2|_P^{k} + S_2|_W^{k}\right) \nonumber \\
+ \Delta tS_2|_P^{k}\left(- a_2 +  c_2S_1|_{P}^{k}\right)
\;\;\;\;\;\;\;\;\;\;\;\;\;
\label{calor2dis}
\end{eqnarray}
where $\sigma_j =\frac{\Delta t \; D_j}{\Delta x^2}$ for $j = 1,2$ .

Note that the equations \ref{calor1dis} and \ref{calor2dis} are non-linear, due to the reactive terms $F_1$ and $F_2$. On the other hand, the Von Neumann superposition procedure is valid only for linear equation systems. Then it is necessary to linearize the equations \ref{calor1dis} and \ref{calor2dis}. 

It is assumed, in the nonlinear terms of \ref{calor1dis} and \ref{calor2dis}, that the variables $S_1$ and $S_2$ are locally positive constants $m_1$ and $m_2$, respectively. From this hypothesis, it follows that the linear version of the equations \ref{calor1dis} and \ref{calor2dis} is given by
\begin{eqnarray}
S_1|_P^{k+1} = S_1|_{P}^{k} + \sigma_1\left(S_1|_E^{k} - 2S_1|_P^{k} + S_1|_W^{k}\right) \nonumber \\
+ \Delta tS_1|_P^{k}\left(a_1 - b_1m_1 - c_1m_2\right)
\label{calor1dislinear}
\\ 
S_2|_P^{k+1} = S_2|_{P}^{k} + \sigma_2\left(S_2|_E^{k} - 2S_2|_P^{k} + S_2|_W^{k}\right) \nonumber \\
+ \Delta tS_2|_P^{k}\left(- a_2 +  c_2m_1\right).
\;\;\;\;\;\;\;\;\;
\label{calor2dislinear}
\end{eqnarray}

Assuming that $\left.S_{1} \right|_{P}^{k} = \left.S_{2} \right|_{P}^{k} = e^{\gamma k}e^ {I\xi P}$ in \ref{calor1dislinear} and \ref{calor2dislinear}, we obtain for the point $P$, in the time $k$,
\begin{eqnarray}
e^{\gamma} = 1 - 4\sigma_1\sin^2 \frac{\xi}{2} + \Delta t\left(a_1 - b_1m_1 - c_1m_2\right)\label{vncalor1} \\
e^{\gamma} = 1 - 4\sigma_2\sin^2 \frac{\xi}{2} + \Delta t\left( - a_2 + c_2m_1\right). \;\;\;\;\;
\label{vncalor2}
\end{eqnarray}
The stability condition requires that $|e^{\gamma}| \leq 1$, so
\begin{eqnarray}
- 1 \leq 1 - 4\sigma_1\sin^2 \frac{\xi}{2} + \Delta t\left(a_1 - b_1m_1 - c_1m_2\right) \leq 1\label{11} \\
- 1 \leq  1 - 4\sigma_2\sin^2 \frac{\xi}{2} + \Delta t\left( - a_2 + c_2m_1\right) \leq 1.\label{12} \;\;\;\;
\end{eqnarray}
Finally, after some algebraic manipulations, we get
\begin{eqnarray}
0 \leq \Delta t \leq \frac{\Delta x^2}{2D_1\sin^2 \frac{\xi}{2} - \frac{\Delta x^2}{2}\left(a_1 - b_1m_1 - c_1m_2\right)}
\label{12a} \\
0 \leq \Delta t \leq \frac{\Delta x^2}{2D_2\sin^2 \frac{\xi}{2} - \frac{\Delta x^2}{2}\left(- a_2 + c_2m_1\right)},
\label{12b} \;\;\;\;
\end{eqnarray}
which are the Von Neumann conditions for the predator-prey discretized equations \ref{calor1dis} and \ref{calor1dis}, respectively. Note that the conditions \ref{12a} and \ref{12b} impose restrictions on the spatial refinements $\Delta x$ and temporal refinements $\Delta t$ for stability in the numerical experiments.

\subsection*{Von Neumann conditions for a Telegraph system}\label{sectionsera}
Consider the telegraph equation system
\begin{eqnarray}
\tau_{j}\frac{\partial^2 S_{j}}{\partial t^2} + \frac{\partial S_{j}}{\partial t}  = D_{j}\frac{\partial^2 S_{j}}{\partial x^2},\label{Von1}
\end{eqnarray}
where $\tau_j>0$ and $D_j>0$ are delay and diffusion constants, respectively, for $j = 1,2$.

Using the discretization schemes given in \ref{eqsjtt}-\ref{eqsjt} in \ref{Von1}, it follows that
\begin{eqnarray}
\left.S_{j} \right|_{P}^{k+1} = \left(2 - \frac{\Delta t}{\tau_j} - \frac{2D_j{\Delta t}^2}{{\Delta x}^2\tau_j}\right)\left.S_{j} \right|_{P}^{k} \;\;\;\;\;\;\;\;\;\;\; \nonumber \\
+ \frac{D_j{\Delta t}^2}{\Delta x^2\tau_j}\left(\left.S_{j} \right|_{E}^{k} + \left.S_{j} \right|_{W}^{k}\right) - \left(1 - \frac{\Delta t}{\tau_j}\right)\left.S_{j} \right|_{P}^{k-1}.\hspace{1mm}\label{Vondis}
\end{eqnarray}
From the von Neumann hypothesis $\left.S_{j} \right|_{P}^{k} = e^{\gamma k}e^{I\xi P}$, we get that 
\begin{eqnarray}
e^{\gamma} = 2 - \frac{\Delta t}{\tau_j} - \frac{2D_j\Delta t^2}{\tau_j\Delta x^2}\left(1 - \cos{\xi} \right) - \left(1 - \frac{\Delta t}{\tau_j}\right)e^{-\gamma}.\label{Von1_1}
\end{eqnarray}
Multiplying \ref{Von1_1} by $e^{\gamma}$ and making use of trigonometric transformations, equation \ref{Von1_1} is rewritten as
\begin{eqnarray}
e^{2\gamma} - 2\beta_j e^{\gamma} + \left(1 - \frac{\Delta t}{\tau_j}\right) = 0,\label{Von1_2}
\end{eqnarray}
where $\beta_j = 1 - 2\sigma_j\frac{\Delta t}{\tau_j}\left(\sin^2{\frac{\xi}{2}}\right) - \frac{\Delta t}{2\tau_j}$ and $\sigma_j =  \frac{D_j\Delta t}{\Delta x^2}$, for $j = 1, 2$. Solving the roots of \ref{Von1_2}
\begin{equation}
    g_{1,2} = e^{\gamma} = \beta_j \pm \sqrt{\beta_j^2 - 1 + \frac{\Delta t}{\tau_j}}.\label{nutela}
\end{equation}

To analyze the roots \ref{nutela}, consider the following lemma.

\vspace{0.5cm}
\noindent
LEMMA I: Let $x \in (0.1]$ and $a(x)$ be a real function defined by $a(x) = \frac{x}{2}$. If for each $x \in (0.1 ]$ we take $y \leq 1 - a(x)$, then
  $0 \leq \sqrt{y^2 - 1 + x} + y \leq 1$.
\vspace{0.5cm}
  
\noindent
Proof: Let $y$ be a real number such that $y \leq 1 -\frac{x}{2}$ for all $x \in (0,1]$. Thus, we have $y \in \left[\frac{1}{2},1\right)$. So,
  \begin{eqnarray*}
  0 \leq\sqrt{y^2 - 1 + x} + y & \leq & \sqrt{ \left(1 - \frac{x}{2}\right)^2 - 1 +x} + 1 -\frac{x}{2} \\
  & = & \sqrt{\frac{x^2}{4}} + 1 - \frac{x}{2}=1.
  \end{eqnarray*}
Therefore, $0 \leq \sqrt{y^2 - 1 + x} + y \leq 1$.
 $\hfill\Box$
\vspace{0,5cm}

In the context of LEMMA I, consider the following cases:

\vspace{0,5cm}
\noindent
Case I) \hspace{0,3cm} $\beta_j^2-1 + \frac{\Delta t}{\tau_j} < 0$.

\noindent
In this case $\sqrt{\beta_j^2 - 1 + \frac{\Delta t}{\tau_j}}$ generates conjugate complex roots, that is,
 \begin{equation}
  g_{1,2} = \beta_j \pm i\sqrt{1 - \frac{\Delta t}{\tau_j} - \beta_j^2}\; .\label{ven}
 \end{equation}
 Taking the norm of\ref{ven}
 \begin{equation}
    ||g||^2 = \beta_j^2 + 1 - \frac{\Delta t}{\tau_j} -\beta_j^2 = 1 - \frac{\Delta t}{\tau_j} < 1 
    \hspace{0,05cm}, \mbox{if} \hspace{0,2cm} 0 \leq \frac{\Delta t}{\tau_j} \leq 1.
  \end{equation}
 So, in this case, when  $0 \leq \frac{\Delta t}{\tau_j} \leq 1$ and $\beta_j^2 - 1 + \frac{\Delta t}{\tau_j} < 0$, by Von Neumann's criterion there is numerical convergence for the \ref{Vondis} system.
  
\vspace{0,5cm}
\noindent
Case II) \hspace{0,3cm} $\beta_j^2-1 + \frac{\Delta t}{\tau_j} \geq 0$. 

\noindent
In this case, for $0 \leq \frac{\Delta t}{\tau_j} \leq 1$, equation (62) can result in $g_{j} = e^{\gamma} = 2\beta_j$, when $\frac{\Delta t}{\tau_j}=1$, therefore we must limit the possible values of $\beta_j$. From LEMMA I, we impose that
\begin{equation}
\frac{1}{2} \leq \beta_j\leq 1 -
\frac{\Delta t}{2\tau_j},
\end{equation}
which results in
\begin{equation}
\frac{1}{2} \leq  1 - 2\sigma_j\frac{\Delta t}{\tau_j}\left(\sin^2{\frac{\xi}{2}}\right) \leq 1.
\end{equation}
Solving the inequality (66), it follows that
\begin{eqnarray}
0 \leq \sigma_j \leq \frac{1}{4}.
\label{condicaoSigma1}
\end{eqnarray}
Then, by LEMMA I, if $0 < \frac{\Delta t}{\tau_j} \leq 1$ and $\beta_j^2-1 + \frac{\Delta t}{\tau_j} \geq 0$, the Von Neumann stability condition of the discretized equation \ref{Vondis} is given by
\begin{eqnarray}
0 \leq \sigma_j \leq \frac{1}{4},
\label{condicaoSigma11}
\end{eqnarray}
where $\sigma_j =  \frac{D_j\Delta t}{\Delta x^2}$ for $j = 1, 2$.

\subsection*{Stability Diagram for a Telegraph Predator-Prey System}\label{secaovixe}
In the case of equations \ref{EQ:R_T_C} and \ref{EQ:R_T_C_1}, or in its discretized form given by the equations \ref{finalpresa} and \ref{finalpredador}, which describe a Telegraph Predator-Prey system with diffusive, reactive, and delayed effects, it was not possible to obtain an analytical form for the Von Neumann stability condition. In this context, numerical experiments will be carried out to obtain the stability diagram for the numerical scheme used.

Setting the reactive parameters $a_1 = 1 \;s^{-1}$, $a_2 = 0.75 \;s^{-1}$, $b_1 = 0.5 \; m/s$, $c_1 = 0.5 \; m/s$, and $c_2 = 0.5 \; m/s$, phase diagrams of the variables $D_j$ (diffusibility coefficient) and $\tau_j$ (delay time) were constructed, for $j = 1, 2$. The phase diagrams show the regions of stability and instability of the our PDEs system, for different discretizations of $\Delta t$ and $\Delta x$. It is important to note that $D_1 = D_2$ and $\tau_1 = \tau_2$ in all simulations, that is, there was no different diffusion or delay between the populations $S_1$ and $S_2$. In the following numerical experiments, the following initial and boundary conditions were considered
 \begin{eqnarray}
      S_{1}(x,0) = S_{1}^{0}=S_{2}(x,0) = S_{2}^{0}=
      \begin{cases}
15 &\mbox{if} \hspace{0.4cm} 24 \leq x \leq 26\\
0 &\mbox{if} \hspace{0.4cm} \mbox{otherwise} \\
      \end{cases}
\label{condicaoinicial1_1}\\ \nonumber  \\
S_{1}(0,t)=S_{1}(50,t)=S_{2}(0,t)=S_{2}(50,t) = 0, 
\;\;\;\;\;\;\;\;\;
\label{condicaodecontorno_1}
       \end{eqnarray}
where $x\in[0,50]$ and $t\in[0,100]$.

All the simulations contained in this work were carried out in an Ubuntu operating system, where the codes were all programmed in FORTRAN 90 language and the images generated via gnuplot software. 

First, some simulations were carried out with different partitions of the temporal domain. Fixed $\Delta x = 0.1 $, Figures 3 and 4 show the graphs of the numerical stability/instability regions, as a function of the diffusion ($D_j$) and the delay time ($\tau_j$) of prey and predator populations, for two different values of $\Delta t$, $\Delta t=0.002$ in Figure 3 and $\Delta t=0.001$ in Figure 4.
\begin{figure}[H]
	\begin{center}
	\caption{Region of stability/instability for the Telegraph Predator-Prey system for refinements $\Delta x = 0.1 $ and $\Delta t = 0.002$.}
\includegraphics[width=9.5cm]{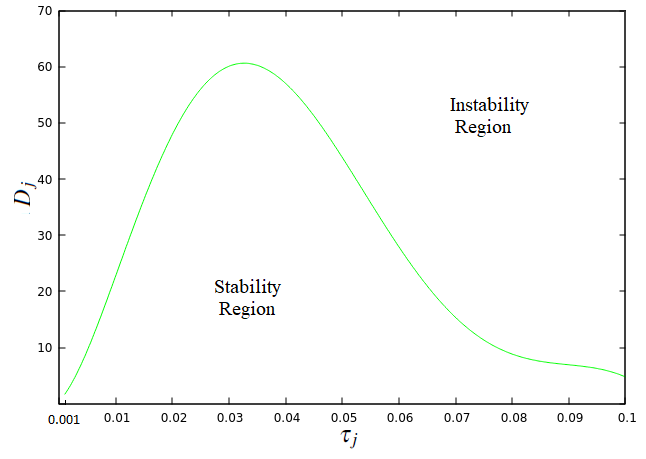}
\label{lol11}
	\end{center}

\end{figure}
\begin{figure}[H]
	\begin{center}
	\caption{Region of stability/instability for the Telegraph Predator-Prey system for refinements $\Delta x = 0.1 $ and $\Delta t = 0.001$.}
\includegraphics[width=9.5cm]{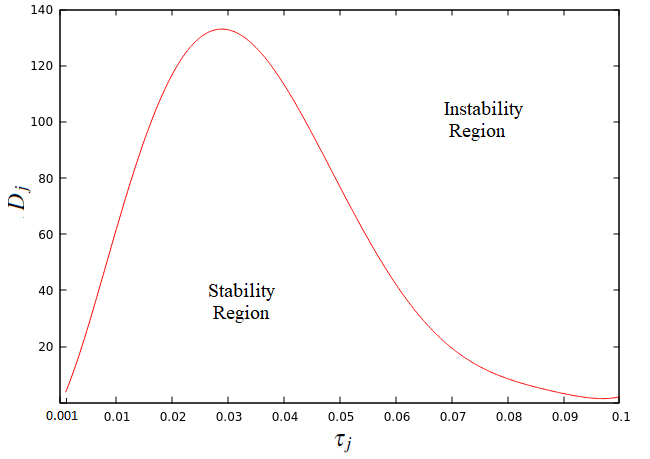}
\label{lol22}
	\end{center}
\end{figure}

Note that when refining the temporal domain, to a fixed $\Delta x$, the numerical stability region grows. See the values of $D_j$ in Figures 3 and 4.

Now we set $\Delta t = 0.0015255$. Figures 5 and 6 show the graphs of the numerical stability/instability regions, as a function of the diffusion ($D_j$) and the delay time ($\tau_j$) of prey and predator populations, for two different values of $\Delta x$, $\Delta x = 0.1$ in Figure 5 and $\Delta x=0.025$ in Figure 6.

\begin{figure}[H]
	\begin{center}
	\caption{Region of stability/instability for the Telegraph Predator-Prey system for refinements $\Delta x = 0.1 $ and $\Delta t = 0.0015255$.}
\includegraphics[width=10.0cm]{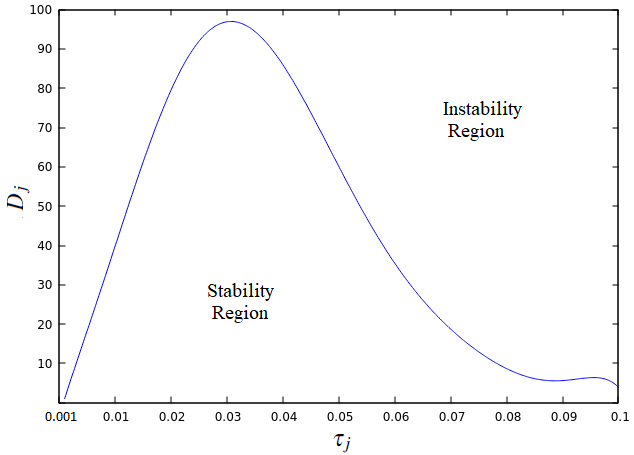}
\label{lol1}
	\end{center}
\end{figure}
\begin{figure}[H]
	\begin{center}
	\caption{Region of stability/instability for the Telegraph Predator-Prey system for refinements $\Delta x = 0.025$ and $\Delta t = 0.0015255$.}
\includegraphics[width=10.0cm]{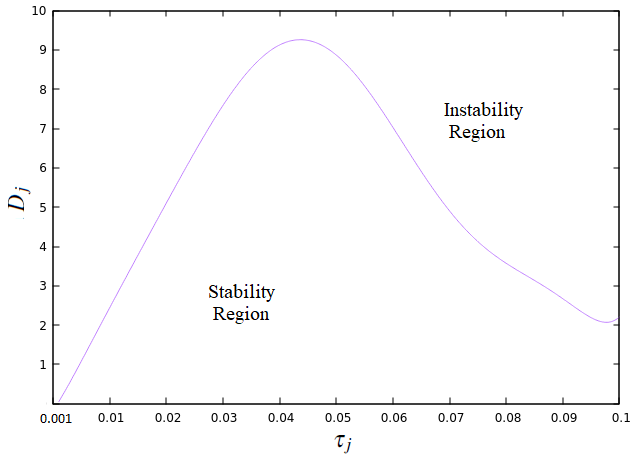}
\label{lol123}
	\end{center}
\end{figure}

Note that when refining the spatial domain, to a fixed $\Delta t$, the numerical stability region decreases. See the values of $D_j$ in Figures 5 and 6.

\section*{Numerical Convergence Analysis}
The objective of this section is to analyze the convergence of the numerical scheme used to solve the model under study. In general, to verify numerical convergence, the Lax Equivalence Theorem must be satisfied \cite{lalax,laxlutor,romeiro1}.

\vspace{0.3cm}
\noindent
Lax Equivalence Theorem:
For a consistent finite difference scheme of a well-posed initial value problem, stability is a necessary and sufficient condition for convergence.
\vspace{0.0cm}

First, it is necessary to verify if the PDEs system are well-posed, that is, if the problem in question has a solution and if it is unique. In the literature there are works on the existence and uniqueness of solutions of Telegraph equations \cite{telegraf}. There are also several works on the existence and uniqueness of Telegraph Predator-Prey systems, or more general Telegraph models  \cite{todabem1,toda,toda1,todabem,toda2}.

Then, a numerical study of mesh refinement is presented to evaluate the process of numerical convergence.

\subsection*{Mesh Refinement and Convergence}
Note that $S_{1}$ and $S_{2}$ are population densities of prey and predators, respectively, so to obtain the total population of prey and predators, in a time $t$, one must integrate $ S_{1}$ and $S_{2}$ at the position coordinate $x$. This calculation is performed using the definite integral
\begin{equation}
P_{j}(t) = \int\limits_{X_{INI}}^{X_{FIN}} S_{j}(x,t)dx,
\end{equation}  
where $P_{1}$ and $P_{2}$ are the populations of prey and predator, respectively, at time $t$.

Let $NI$ be the number of partitions in the spatial domain and $NJ$ the number of partitions in the temporal domain. Thus, we have the following relations
\begin{eqnarray}
\Delta t = \frac{T_{FIN} - T_{INI}}{NJ - 1}\label{dt}
\;\;\;\; \mbox{and} \;\;\;\;
\Delta x = \frac{X_{FIN} - X_{INI}}{NI - 1}.\label{dx}
\end{eqnarray}
Next, an analysis of the convergence of the numerical solution of the problem regarding the mesh refinement is performed. The Table \ref{tabelatestetau0} presents the values of the parameters used in the simulations in this subsection. The initial and boundary conditions are those given earlier in \ref{condicaoinicial1_1} and \ref{condicaodecontorno_1}. 
\begin{table}[H]
\renewcommand{\arraystretch}{1.2}
\centering
\caption{Values for the model parameters used in the numerical simulations presented in Table 2.}
\label{tabelatestetau0}
\begin{tabular}{ccc}
\hline
\textbf{Parameters}                                      & \textbf{Prey } & \textbf{Predator } \\ \hline
\textbf{$a_{j} \; (s^{-1})$}       & 1.0                    & 0.75                        \\ 
\textbf{$b_{j} \; (m/s)$}       & 0.5                    & 0.0                        \\ 
\textbf{$c_{j} \; (m/s)$}        & 0.5                    & 0.5                        \\ 
\textbf{$D_{j} \; (m^2/s)$}         & 1.0                    & 1.0                        \\
\textbf{$\tau_{j} \; (s)$}         & 0.001                    & 0.001                        \\ 
$ T_{INI} \; (s)$              & 0.0                  & 0.0                      \\ 
$ T_{FIN} \; (s)$              & 100.0                  & 100.0                      \\
$ X_{INI} \; (m)$              & 0.0                  & 0.0                      \\ 
$ X_{FIN} \; (m)$              & 50.0                  & 50.0 \\ \hline
\end{tabular} \\
\end{table}

Table \ref{tabelatudo} presents numerical experiments for the solution of the discretized Telegraph Predator-Prey system (23-32), as a function of the mesh refinement. The process of population convergence is observed as a function of this refinement.
\begin{table}[H]
\renewcommand{\arraystretch}{1.2}
\centering
\caption{Prey population $P_{1}(t) $ and predator population $P_{2}(t) $  of the Telegraph Predator-Prey model for various refinements $\Delta x$ and $\Delta t$, at time $t = 100 \; s$, with the parameters in Table 1.}
\label{tabelatudo}
\begin{tabular}{cccc}
\hline
\textbf{$\Delta x$ $(m)$}                                      & \textbf{$\Delta t$ $(s)$ } & \textbf{$P_{1}(t=100) $}   & \textbf{$P_{2}(t=100) $}\\ \hline
5.0    & 0.002               & 62.69325     & 18.90470        
\\ 
3.5    & 0.002               & 65.13776     & 19.45726        
\\ 
2.0    & 0.002               & 70.27077     & 20.51610        
\\ 
1.0    & 0.002               & 71.68832     & 21.04255        
\\ 
0.75    & 0.00175               & 72.22790     & 21.17758         
\\ 
0.5    & 0.0015               & 72.56129    & 21.29889          
\\
0.1    & 0.001               & 73.13264    & 21.49681           
\\ 
 0.05       & 0.0005              & 73.20856   & 21.52185                
\\ 
0.04        & 0.0004              & 73.22319  & 21.52683              
\\ 
0.03         & 0.0003          & 73.23832       & 21.53186              
\\
0.02       & 0.0002              & 73.25332    & 21.53683                        \\ 
0.01        & 0.0001              & 73.26829    & 21.54183                        \\ 
0.005        & 0.00005              & 73.27579    & 21.54433                        \\ 
0.0025        & 0.000025              & 73.27954    & 21.54558                        \\ \hline
\end{tabular}\\
\end{table}

\subsection*{Numerical Experiments on Convergence in the Telegraph Predator-Prey System}\label{simu1}
In this subsection, two numerical experiments are carried out
for the Telegraph Predator-Prey system (11-12) in order to better understand the stability/instability phase diagrams obtained in Figures 3-6. 

Let us consider the stability and instability regions obtained in Figure 5, when $\Delta t=0.0015255$ and $\Delta x=0.1$. Figure 7 presents these regions again. Figure 7 also shows two points, one in the region of stability of the numerical scheme and the other at the limit of the region of instability of the numerical scheme.
\begin{figure}[h]
    \begin{center}
    \caption{Region of stability/instability for the Telegraph Predator-Prey system for refinements $\Delta x = 0.1$ and $\Delta t = 0.0015255$. The point $ (0.05,20)$ is in the stability region, while the point $ (0.05,61)$ is at the limit
of the instability region.}
\includegraphics[width=10.0cm]{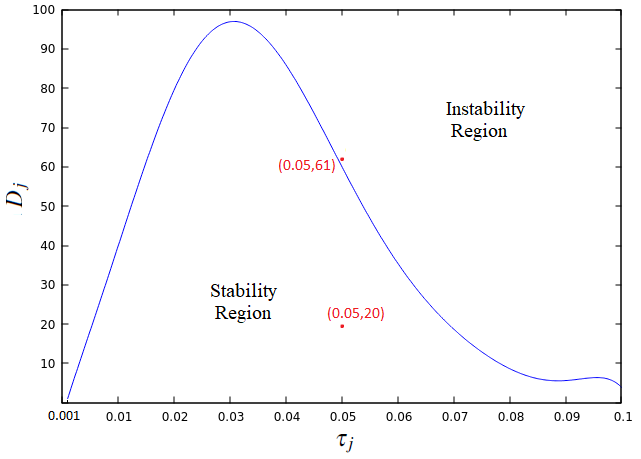}
\label{oae}
	\end{center}
\end{figure}

For the first numerical experimentation in this subsection, the values of the model parameters are presented in Table \ref{tabelatestetau111}. Note that this simulation corresponds to the point in the stability region of Figure 7.  
\begin{table}[H]
\renewcommand{\arraystretch}{1.2}
\centering
\caption{Model parameter values for the first numerical experiment, corresponding to the point in the stability region of Figure 7.}
\label{tabelatestetau111}
\begin{tabular}{ccc}
\hline
\textbf{Parameters}                                      & \textbf{Prey }  & \textbf{Predator }  \\ \hline
\textbf{$a_{j} \; (s^{-1})$}       & 1.0                    & 0.75                        \\ 
\textbf{$b_{j} \; (m/s)$}       & 0.5                    & 0.0                        \\ 
\textbf{$c_{j} \; (m/s)$}        & 0.5                    & 0.5                        \\ 
\textbf{$D_{j} \; (m^2/s)$}         & 20.0                    & 20.0                        \\
\textbf{$\tau_{j} \; (s)$}        & 0.05                    & 0.05                        \\ \hline
\end{tabular} \\
\end{table}

\begin{figure}[ht]
\begin{center}
\caption{Population densities of prey and predators of the Telegraph Predator-Prey model, in $t=100$, corresponding to the first numerical experiment.}
\includegraphics[width=10.0cm]{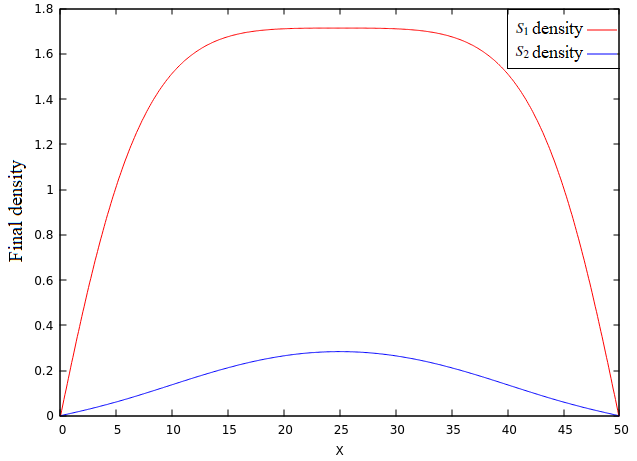}
\label{log2a}
\end{center}
\end{figure}

Figure \ref{log2a} presents the population densities of prey and predators, in space coordinate, for the first numerical experiment at time $t=100$.

\begin{figure}[H]
\begin{center}
\caption{Populations of prey and predators over time of the Telegraph Predator-Prey model, in $x=25$, corresponding to the first numerical experiment.}
\includegraphics[width=9.5cm]{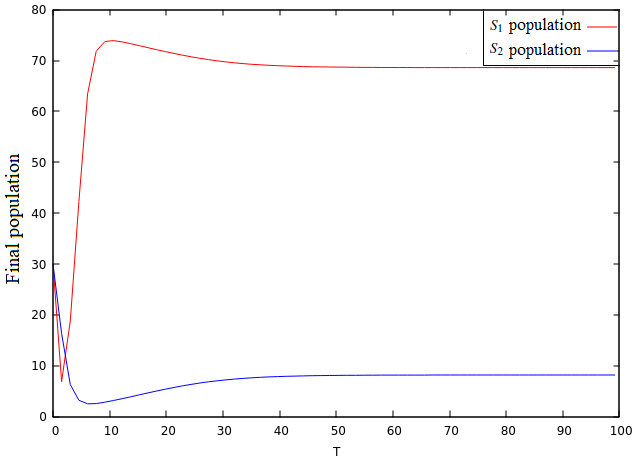}
\label{log2b}
\end{center}
\end{figure}

Figure \ref{log2b} presents the temporal variation, in $x=25$, of the populations of prey and predators in the first numerical experiment.

\begin{figure}[H]
\begin{center}
\caption{3D graph of predator and prey population densities of the Telegraph Predator-Prey model corresponding to the first numerical experiment.}
\includegraphics[width=11.0cm]{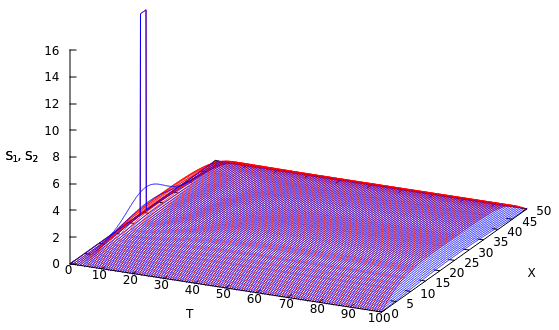}
\label{log22}
\end{center}
\end{figure}

Finally, Figure \ref{log22} shows population densities of predator and prey in the entire spatial and temporal domain. Note that in this simulation there are no negative population densities.

For the second numerical experiment in this subsection, the values of the model parameters are presented in Table \ref{tabelatestetau233}. The second numerical experiment corresponds to the point at the limit of the instability region in Figure 7.
\begin{table}[H]
\renewcommand{\arraystretch}{1.2}
\centering
\caption{Model parameter values for the second numerical experiment, corresponding to the point at the limit of the instability region in Figure 7.}
\label{tabelatestetau233}
\begin{tabular}{ccc}
\hline
\textbf{Parameters}                                      & \textbf{Prey } & \textbf{Predator } \\ \hline
\textbf{$a_{i} \; (s^{-1})$}       & 1.0                    & 0.75                        \\ 
\textbf{$b_{i} \; (m/s)$}       & 0.5                    & 0.0                        \\ 
\textbf{$c_{i} \; (m/s)$}        & 0.5                    & 0.5                        \\
\textbf{$D_{i} \; (m^2/s)$}         & 61.0                    & 61.0                        \\ 
\textbf{$\tau_{i} \; (s)$}        & 0.05                    & 0.05                        \\ \hline
\end{tabular} \\
\end{table}

\begin{figure}[ht]
\caption{Population densities of prey and predators of the Telegraph Predator-Prey model, in $t=100$, corresponding to the second numerical experiment.}
\begin{center}
\includegraphics[width=10.0cm]{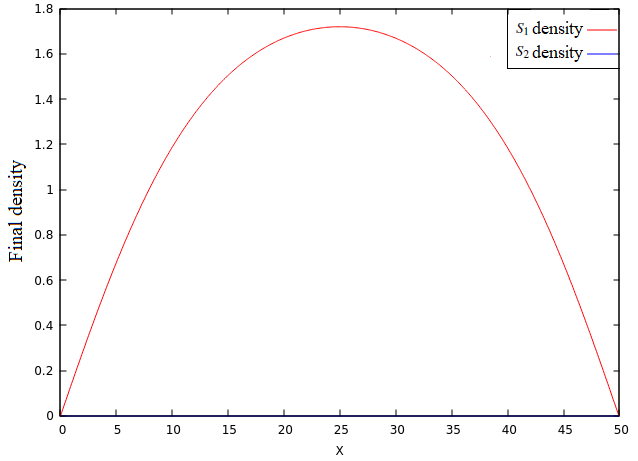}
\label{fig:gull}
\end{center}
\end{figure}

Figure \ref{fig:gull} presents the population densities of prey and predators, in space coordinate, for the second numerical experiment at time $t = 100$. Note that the final population density of predators is zero, that is, the extinction of predators occurs. 

\begin{figure}[H]
\begin{center}
\caption{Populations of prey and predators over time of
the Telegraph Predator-Prey model, in x = 25 , correspond-
ing to the second numerical experiment.}
\includegraphics[width=10.0cm]{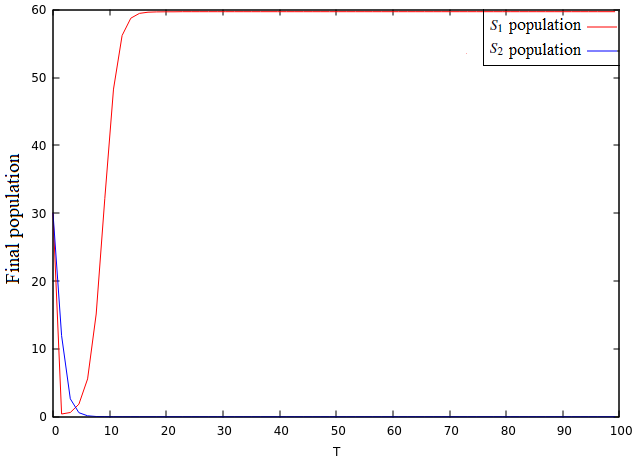}
\label{fig:tiger}
\end{center}
\end{figure}

Figure \ref{fig:tiger} presents the temporal variation, in x = 25, of prey and predators populations in the second  numerical experiment. 

\begin{figure}[H]
\begin{center}
\caption{3D graph of predator and prey population densities of the Telegraph Predator-Prey model corresponding to the second numerical experiment.}
\includegraphics[width=11.2cm]{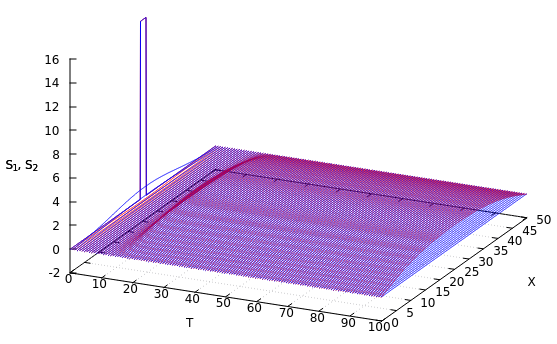}
\label{fig:tiger2}
\end{center}
\end{figure}

\begin{figure}[ht]
\begin{center}
\caption{Observation (first red line) of the spatial/temporal region where negative solutions for prey population density occur.}
\includegraphics[width=10.0cm]{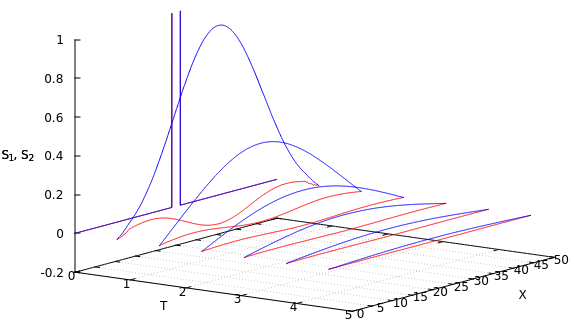}
\label{zoom}
\end{center}
\end{figure}

The Figure \ref{fig:tiger2} shows population densities of predator and prey in the entire spatial and temporal domain. It is observed that for some values of $t$ there are negative population densities, as can be seen in Figure \ref{zoom}.

Note this behavior is unrealistic from a biological point of view. According to the numerical solution presented in the second experiment, given the initial condition (69), the prey population rapidly decreases, becomes negative and grows again until it reaches saturation. This is evidence of the instability of numerical solutions.

\section*{Conclusion}
In this work, the mathematical modeling of the   diffusive-reactive predator-prey system with delay (11)-(12) was performed. This model is also called Telegraph Predator-Prey system. Through numerical modeling, this EDP system was discretized by the Finite Difference method, equations (23)-(32).

First, it was verified that the discretized equations (23)-(32) were consistent with the PDEs (11)-(12) of the mathematical model.

In the sequence, through the Von Neumann procedure, the numerical stability of systems simpler than the Telegraph Predator-Prey system was discussed. For the predator-prey system (46)-(47) and for the telegraph equation system (58) it was found that the numerical scheme used is conditionally stable, since the stability constraints (56)-(57) and (68) depend on the parameters of these models. 

On the other hand, for the telegraph predator-prey system (11)-(12), it was not possible to obtain an analytical mathematical relationship for the Von Neumann stability condition. Then, through numerical experiments, phase diagrams of the stability/instability regions were constructed as a function of the diffusion $D_j$ and delay time $\tau_j$ parameters of the model (11)-(12), see Figures 3-6. Thus, it can be concluded that the telegraph predator-prey system (11)-(12) is also conditionally stable, depending on its diffusion, delay and reactive parameters.

Through numerical experiments again, the convergence of the numerical model was analyzed. Through the spatial and temporal refinement of the mesh, it was found that the prey and predator populations of the telegraph predator-prey system converged, as shown in Table 2.

Finally, it was possible to observe the numerical instabilities that occur in numerical simulations. In the phase diagram presented in Figure 7, two points were chosen, one in the region of stability of the numerical scheme and the other at the limit of the region of instability of the numerical scheme. It was found that when passing from the stability region to the instability region, negative fluctuations occur for population densities, which is an evidence of instability of numerical solutions, as can be seen in Figure 14.

\section*{Acknowledgments}
The authors Kariston Stevan Luiz and Juniormar Organista was partly financed by the Coordination of Improvement of Higher Education Personnel (Capes), Brazil - Finance Code 001.

{}

\end{document}